\documentclass[letterpaper, 10 pt, conference]{ieeeconf}  

\IEEEoverridecommandlockouts                              

\title{\LARGE \bf
Diagnosis of Anomaly in the Dynamic State Estimator of a Power System using System Decomposition
}

\author{Malini Ghosal\\
Electrical and Computer Engineering\\
Texas Tech University, Lubbock, TX 79409 USA\\
Email: {\tt malini.ghosal@ttu.edu}
\thanks{*The material is based upon the work performed by the author at Texas Tech University, under the partial support of National Science Foundation Grant No. DGE-1438921.}
}

\usepackage{color}

\usepackage{amsmath,amssymb,graphicx,graphics}
\usepackage{bm}
\usepackage{subfigure}
\newtheorem{proposition}{Proposition}
\newtheorem{definition}{Definition}

\newtheorem{example}{Example}

\newtheorem{assumption}{Assumption}
\usepackage{algorithm,algorithmic}
\usepackage{ragged2e}

\def\0{\mathbf{0}}

\def\x{\mathbf{x}}

\def\O{\mathbf{\mathcal{O}}}

\def\xbigbig{\mathbf{\mathcal{X}}}

\usepackage{listings}

\usepackage{cite}
\begin{document}

\maketitle
\thispagestyle{empty}
\pagestyle{empty}

\begin{abstract}
In a state estimator, the presence of malicious or simply corrupt  sensor data or bad data is detected by the high value of normalized measurement residuals that exceeds the threshold value, determined by the $\chi^2$ distribution. However,  high normalized residuals can also be caused by another type of anomaly, namely gross modeling or topology error. In this paper we propose a method to distinguish between these two sources of anomalies - 1) malicious sensor data and 2) modeling error. The anomaly detector will start with assuming a case of malicious data and suspect some of the individual measurements corresponding to the highest normalized residuals to be `malicious', unless proved otherwise. Then, choosing a change of basis, the state space is transformed and decomposed into `observable' and `unobservable' parts with respect to these `suspicious' measurements. We argue that, while the anomaly due to malicious data can only affect the `observable'  part of the states, there exists no such restriction for anomalies due to modeling error. Numerical results illustrate how the proposed anomaly diagnosis based on Kalman decomposition can successfully distinguish between the two types of anomalies.   
\end{abstract}


\section{Introduction}
With the high frequency synchronized measurements from the Phasor Measurement Units (PMUs), the real-time and fast changing dynamic states  of an interconnected power system can be estimated using Kalman filtering based techniques\cite{zhou2014capturing, Zhang:2011, feasibility, Bila:2014}. Starting from a known operating condition, the (noisy) dynamic states at each time-step are estimated by optimally combining the dynamical model-based predicted state-information with the real-time measurement data \cite{Simon:2006}.
Therefore the presence of bad measurement data results in  highly inaccurate  state information; subsequently affecting the control and operating decisions. Such bad data are  detected by inspecting the normalized measurement residual and performing the $\chi^2$ hypothesis tests \cite{Abur:2004,Giani:2013,Mo:2010}. However, these bad data detectors can be  falsely triggered if there are {major configuration changes} due to topology error, significant parameter estimation error or unanticipated large load change etc \cite{Abur:2004, bretas2013geometrical}. It is therefore crucial to distinguish between the underlying causes of anomaly: i.e., the presence of malicious measurement data (scenario-1) or significant modeling error (scenario-2). 

In the power systems literature,  detection of modeling error for special cases has been addressed \cite{lourencco2015topology, weng2015convexification, lourencco2008unified, gou2008new}. However in this {paper} a generic method to {differentiate} between the sources of anomalies, namely bad measurement data and modeling error, is presented, based on Kalman {decomposition} techniques. Topology identification problem has been separately addressed extensively in the literature in \cite{clements1988detection, costa1993identification}.   Separately topology error detection using static state estimation has been addressed in \cite{Abur:2004, lourencco2008unified, de2002implicitly, lourencco2006topology, zhao2017enhanced, weng2015convexification,bretas2013geometrical}. Distinguishing the source of anomaly have been addressed and solved for static state estimation in \cite{lourencco2008unified} using Lagrangian multiplier. In \cite{lourencco2015topology}, the authors detected topology error as a cyber-attack for the dynamic state estimation.

However, distinguishing the source of anomaly for the dynamical state estimation of power system has been understudied. The rest of the paper is organized as below. Section II describes the dynamic state estimation process and the bad data detection method using $\chi^2$-tests on the measurement residuals. Section III describes the problem formulation and Section IV describes the proposed methodology. The algorithm is presented in a concise form in Section\,V\,. Section\,VI presents the simulation results, while the conclusions are drawn in Section\,VII.

\section{Dynamic state estimation and bad data detector}
The power systems is a complex nonlinear dynamical network with multiple time-scales. Typically, under the fast time-scales of operation, e.g. transient dynamics, the power systems can be modeled as a set of nonlinear differential equations. Details on developing such models can be found in \cite{ghosal2017fusion, Bila:2014, Zhang:2011, Kundur:1994}. In this and the following two sections, the problem and the proposed solution would be presented in a generic set-up, applicable for any (linear) dynamical systems. Henceforth, noting that the presented methodology is generic and not specific to power systems, we have omitted the modeling details of power system dynamics for brevity. 

Let us consider a dynamical system described by the following set of nonlinear state equations,
\begin{align}
    \x(k+1) & = \mathbf{g} \lbrace \x(k) \rbrace  \rbrace + \mathbf{w}(k), 
\end{align}
where, $\x \in \mathbb{R}^{n}$ is the state vector, $ \mathbf{g} :\mathbb{R}^n\mapsto\mathbb{R}^n$ is nonlinear function describing the dynamics and $\mathbf{w}\in\mathbb{R}^n$ is the zero-mean process noise vector with covariance $\mathbf{Q}$. The observer or the measurements can be described by the following (nonlinear) equations,
\begin{align}
    \mathbf{z}(k+1)= \mathbf{h} \lbrace \x(k+1) \rbrace + \mathbf{v}(k+1),
\end{align}
where $z\in\mathbb{R}^m$ is the measurement vector and $h:\mathbb{R}^n\mapsto\mathbb{R}^m$ is the nonlinear measurement function. Let us assume, without a loss of generality\footnote{The equilibrium can be always moved to the origin with an appropriate shifting of the state variables.}, that $x=0$ is a stable equilibrium point of the system, such that, under normal operating conditions, the states are all close to the origin. Therefore, under the assumption of (close to) normal operations, the system model can be simplified as a linear time-invariant dynamical system as follows,
\begin{subequations}\begin{align}
\x(k+1) &= \mathbf{A}\x(k)  + \mathbf{w}(k),\\
    \mathbf{z}(k+1) &= \mathbf{H}\x(k+1) + \mathbf{v}(k+1),
\end{align}\end{subequations}
where $\mathbf{A} = \frac{ \partial{\mathbf{g}}}{\partial{\x}} \Huge{|}_{x=0} \in \mathbb{R}^{n \times n}$ is the system Jacobian matrix, and $\mathbf{H} = \frac{ \partial{\mathbf{h}}}{\partial{\x}}\Huge{|}_{x=0} \in \mathbb{R}^{m \times n}$ is the measurement Jacobian matrix.


To estimate the states of this system, the estimator is initialized first by solving a set of algebraic power-flow equations. Then the following Kalman filtering equations are used: 
\begin{subequations}\begin{align}
\textit{Initialization}&: \notag \\
\hat{\x}(0|0) & = \mathbb{E} \lbrace \mathbf{x}(0) \rbrace , \\ 
\mathbf{P}(0|0) & = \mathbb{E} \lbrace \mathbf{x}(0) \!-\! \hat{\x}(0|0) \rbrace  \lbrace \mathbf{x}(0) \!-\! \hat{\x}(0|0) \rbrace ^{T}\!\! \\
\textit{Prediction steps} & : \notag  \\
\hat{\x}(k+1|k) &= \mathbf{A} \hat{\x}(k|k), \\
\mathbf{P}(k+1|k) & =\mathbf{A} \mathbf{P}(k|k)\mathbf{A}^{T} + \mathbf{Q}, \\
\textit{Correction steps} & :  \notag \\
\tilde{\mathbf{z}}(k+1) & = \mathbf{z}(k+1) - \mathbf{H} \hat{\x}(k+1|k), \\
\mathbf{S}(k+1) & = \mathbf{H} \mathbf{P}(k+1|k)\mathbf{H}^{T} + \mathbf{R}, \\
\mathbf{K}(k+1) &= \mathbf{P}(k+1|k)\mathbf{H}^{T}\mathbf{S}^{-1}(k+1),\\
\hat{\x}(k+1|k+1)&= \hat{\x}(k+1|k)  + \mathbf{K}(k+1)\tilde{\mathbf{z}}(k+1), \\
\mathbf{P}(k+1|k+1) & =\left (\mathbf{I} - \mathbf{K}(k+1)\mathbf{H} \right )\mathbf{P}(k+1|k),
\end{align}\end{subequations}
where, the variables are explained below, \\
$\mathbf{\hat{x}} $: estimated state, \\
$\tilde{\mathbf{z}}$: measurement residual, \\
$\mathbf{S}$: covariance of measurement residual,\\
$\mathbf{P} $: the covariance of the estimation error, and \\
$\mathbf{K}$: the optimal Kalman gain.\\

Typically, the authenticity of a measurement data is detected by computing the following expression,
\begin{align}
    c(k+1) &=\tilde{\mathbf{z}}^{T}(k+1)\mathbf{S}^{-1}(k+1)\tilde{\mathbf{z}}(k+1),
\end{align}

$c(k+1)$  is a random variable with a ${\chi}^2$ distribution with degrees of freedom ($m_i\leq m$) equal to the number of independent measurements. {If we choose a confidence value of $p$, then the threshold $TH_{\chi}$} may simply be found using \eqref{E:chi_t}.
The value of $TH_{\chi}$ could be determined by ensuring that under normal conditions (\textit{no attack}) the value of the normalized measurement residuals ($c$) is less than $TH_{\chi}$ with certain high probability (referred to as the \textit{confidence level}). We define the probability of the normalized measurement residual being less than the threshold as,
\begin{align}
    F(TH_{\chi}|m_i)&:=\mathbf{Pr}\lbrace c\leq TH_{\chi}\,\vert\,m_i\rbrace\notag \\ &= \int_o^{TH_{\chi}} \frac{t^{\frac{m_i-2}{2}}\exp\left({-t}/{2}\right)}{{2^{\frac{m_i}{2}} \mathcal{\gamma}}(m_i/2)}\,dt\,, \label{E:Thx1}
\end{align}
where $\gamma ( \cdot ) $ is the Gamma function. Note that $F(\cdot)$ is an increasing function of the threshold value. The threshold value is chosen such that $F(TH_{\chi}|m_i)\geq p$\,. Given a (sufficiently high) confidence value of $p\in(0,1)$, the minimum required value of the threshold $TH_{\chi}$ may be found\footnote{MATLAB command `chi2inv' is used.} using the following equation \cite{MATLAB:2015},
\begin{align}  \label{E:chi_t}
TH_{\chi} = F^{-1}(p|m_i):= \left\lbrace x : F(x|m_i)=p \right\rbrace\,. 
 \end{align}

\begin{definition}

Suppose the measurement data $\mathbf{z}(k+1)$ is manipulated by a malicious agent and a new value $\mathbf{z}^{'}(k+1)$ is placed in an attempt to mislead the estimation. {The malicious measurement data is modeled by},
\begin{align}
    \mathbf{z}^{'}(k+1)=\mathbf{z}(k+1) + \mathbf{a}(k+1),
\end{align} where, $\mathbf{a}(k+1) \in \mathbb{R}^{m}$ is the vector of values added to the original measurement.
\end{definition}
{Then}, the changed equations of the estimator at the (k+1)-th time instant will be,
\begin{align}
    \tilde{\mathbf{z}}^{'}(k+1) & = \mathbf{z}^{'}(k+1) - \mathbf{H} \hat{\x}(k+1|k),\\
    &= \tilde{\mathbf{z}}(k+1)+\mathbf{a}(k+1)\notag
\end{align}

{Therefore, the expression for malicious data detector (or, the normalized measurement residual) is}:

\begin{align}
    c^{'}(k+1) &= \tilde{\mathbf{z}}^{'T}(k+1) \mathbf{S}^{-1}(k+1) \tilde{\mathbf{z}}^{'}(k+1) \notag  
    \\ 
    &= c(k+1)+ \mathbf{a}^{T}(k+1)\mathbf{S}^{-1}(k+1)\mathbf{a}(k+1)\notag\\
    &~ + 2\,\mathbf{a}^{T}(k+1)\mathbf{S}^{-1}(k+1)\tilde{\mathbf{z}}(k+1)
\end{align}
Occurrence of malicious measurement data is successfully detected if $c^{'}(k+1) > TH_{\chi} $.

\section{Problem formulation}
While presence of bad measurement data would result in a large $c^{'}(k+1)$, it can also be caused by modeling error, e.g. when the assumed values of $\mathbf{A}$ and $\mathbf{H}$ deviate largely from those of the real system. This could be caused by various reasons. Some of the reasons of this gross error  are, 
\begin{itemize}
\item topology error: suppose, there is a sudden and unanticipated change in topology of which the topology processor is not yet updated. This will affect both $\mathbf{A}$ and $\mathbf{H}$ {matrices}.
\item change of operating conditions: change in operating condition will move the {initial condition} away from the origin (assumed to the stable operating point). This will eventually affect the system parameters $\mathbf{A}$ and $\mathbf{H}$. If the state estimator is unaware of this updated operating point, then the computed Jacobians $\mathbf{A}$ and $\mathbf{H}$ would not represent the true system.
\item some machine or line parameters, such as damping coefficients of generators, line impedance might change. We should remember, however, that small inaccuracies and parameter uncertainties are already taken care of and modeled as an additive Gaussian process noise $\mathbf{w}(k+1)$.
\end{itemize}

\begin{subequations}\begin{align}
\textit{Prediction steps} & : \notag  \\
\hat{\x}'(k+1|k) &= \mathbf{A}' \hat{\x} (k|k) , \\
\mathbf{P}'(k+1|k) & =\mathbf{A}'\mathbf{P}(k|k)\mathbf{A}'^{T} + \mathbf{Q}, \\
\textit{Correction steps} & :  \notag \\
\tilde{\mathbf{z}}'(k+1) & = \mathbf{z}(k+1) - \mathbf{H}' \hat{\x}'(k+1|k), \\
\mathbf{S}'(k+1) & = \mathbf{H}' \mathbf{P}'(k+1|k)\mathbf{H}'^{\, -1} + \mathbf{R}, \\
\mathbf{K}'(k+1) &= \mathbf{P}'(k+1|k)\mathbf{H}'^{T}\mathbf{S}'^{-1}(k+1),\\
\hat{\x}'(k+1|k+1)&= \hat{\x}'(k+1|k)  + \mathbf{K}'(k+1)\tilde{\mathbf{z}}'(k+1), \\
\mathbf{P}'(k+1|k+1) & =\left( \mathbf{I} - \mathbf{K}'(k+1)\mathbf{H}' \right) \mathbf{P}'(k+1|k).
\end{align}\end{subequations}
Now the malicious data indicator,
\begin{align}
    c'(k+1) &=\tilde{\mathbf{z}}'^{T}(k+1) \lbrace \mathbf{S}'(k+1) \rbrace ^{-1} \tilde{\mathbf{z}}'(k+1),
\end{align}
%
%
where, the prime notation (${'}$) represents the corrupted variables  after one time step after the modeling error. Therefore, every occasion when $c'(k+1)$ exceeds the threshold value due to such modeling errors, false alarms will be raised. These false alarms are harmful as these will prompt the system operators to act in certain ways that could be undesired and therefore, should be avoided. In this section, the underlying cause of $c'(k+1)$ exceeding the threshold $TH_{\chi}$ is detected. 

\begin{figure}[tbh]
 \centering
	\includegraphics[scale=.27]{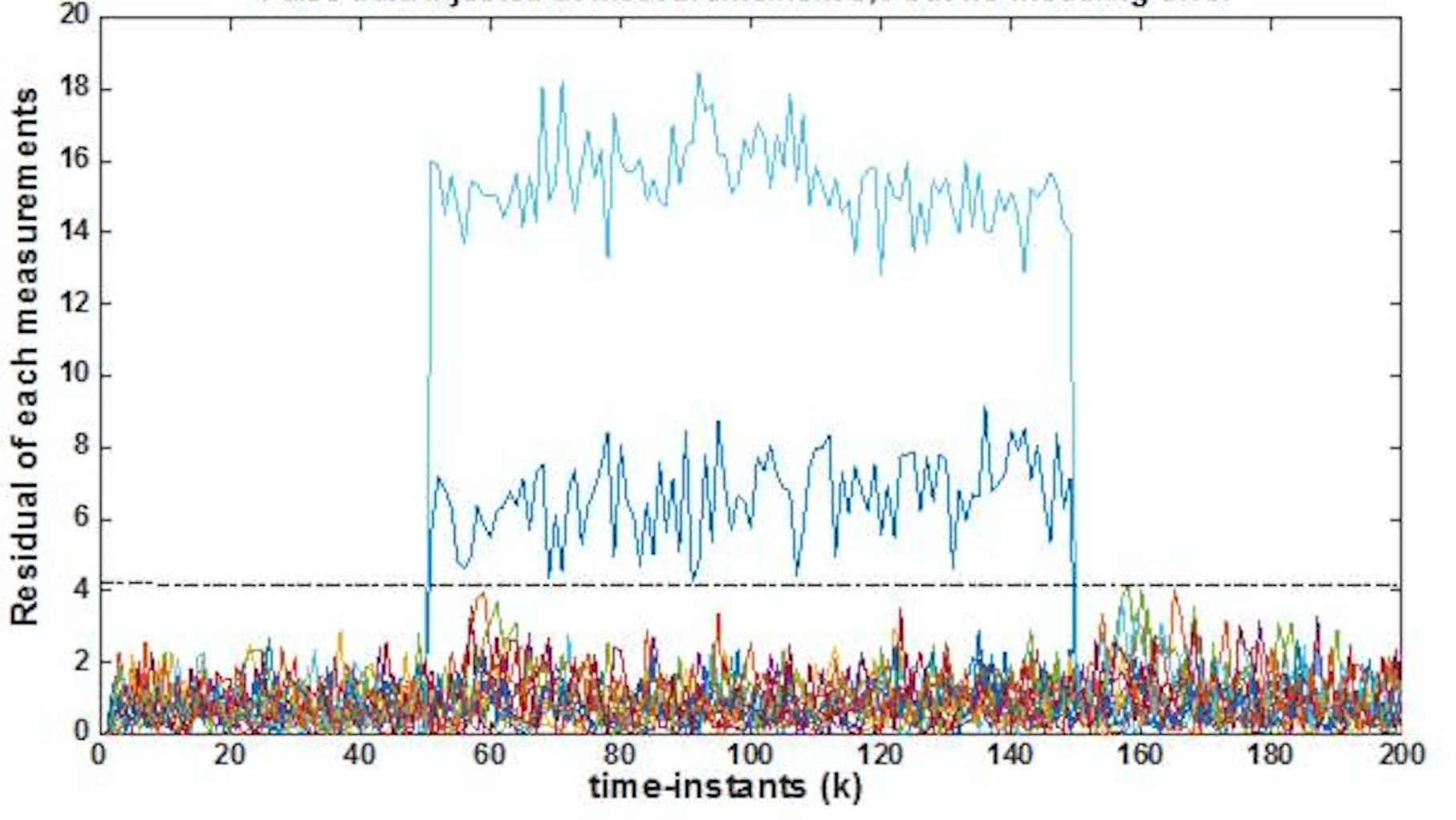}
 \caption{Large residuals due to malicious data. }
 \label{F:AD1}
 \end{figure}

\begin{figure}[tbh]
 \centering
	\includegraphics[scale=0.6]{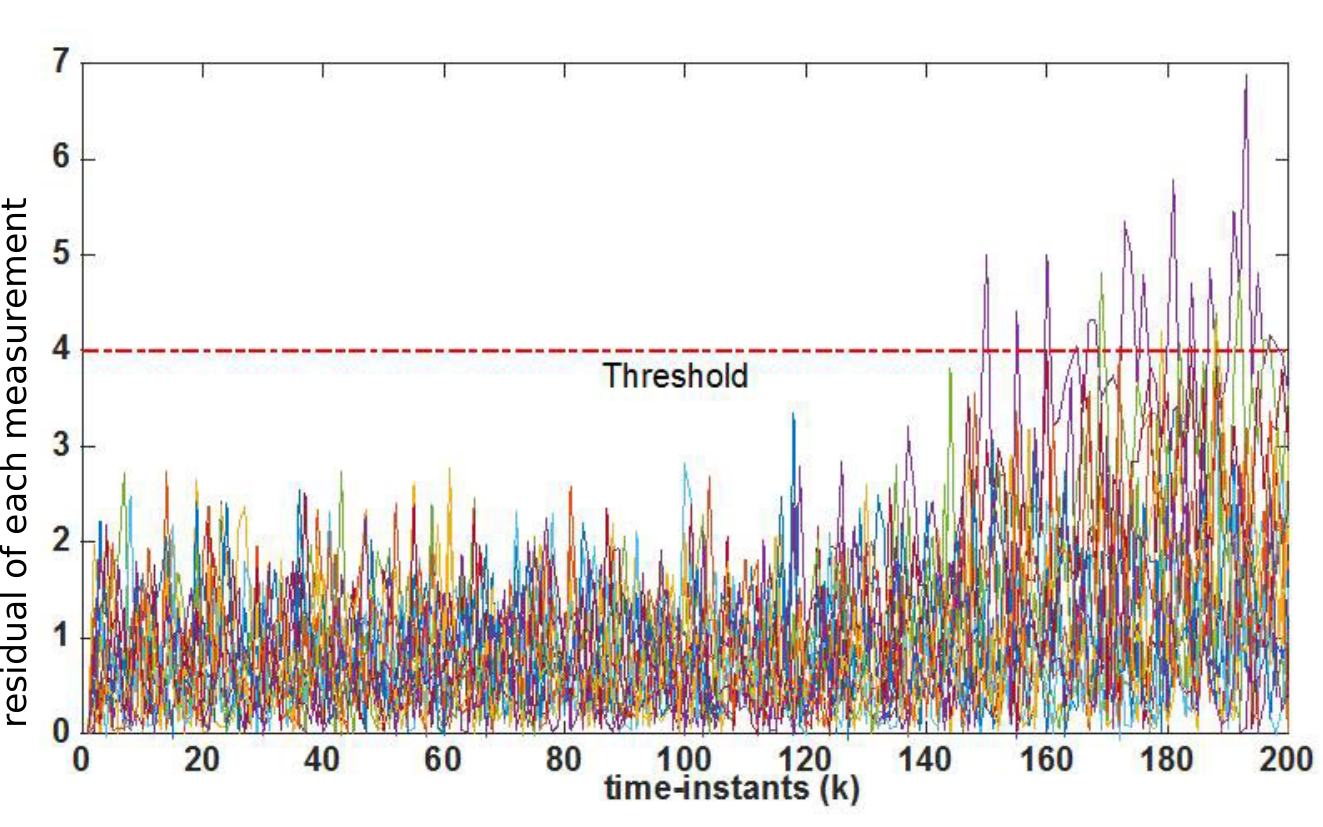}
 \caption{Large residuals due to topology error. }
 \label{F:AD2}
 \end{figure}

\begin{example}
In order to illustrate the anomalies due to bad data and modeling error, an arbitrary tenth order system is simulated. First, malicious data is injected in measurements-3 and 8 during time instants 50-150. From Fig. \ref{F:AD1}, it can be observed that the residuals corresponding to these measurements cross the detection threshold during the attack instants (Scenario-1). The attack vectors injected here have a standard deviation 5 and 10 times  larger than their respective noise variance. Next, in Fig. \ref{F:AD2} the simulation is repeated again. However, instead of any attack injection, a modeling error is simulated: topology is not updated due to a fault occurring after time instant 150. The residuals corresponding to some measurements misleadingly exceed the detection threshold (Scenario-2). Therefore, simply by looking at the normalized residuals, it is not possible to detect the anomaly.
\end{example} 
In the following sections, a method will be introduced to distinguish between such anomalies.

\section{Proposed Methodology}
By inspecting the normalized residuals for all the available measurements at time $k$, denoted by $\mathbf{z}(k)$, suspicious measurements ($\mathbf{z}_s\subset \mathbf{z}$) are identified. Let us assume that only a small subset of the measurements have high residuals such that the observability matrix with respect to the \textit{suspicious measurements} ($\mathbf{z}_s(k) \in \mathbb{R}^{m_s}$) is rank deficient\footnote{We assume that the chances of a large number of measurements being attacked/corrupt concurrently are thin.}. 
Associated with $\mathbf{z}_s$, there exists an unobservable sub-space, denoted by $\mathcal{N}(\mathcal{\mathbf{O}}_s)$\,, and an observable sub-space, denoted by $\mathcal{N}(\mathcal{\mathbf{O}}^{\perp}_s)$\,.
Therefore, an assumption {for the observability matrix shown below is made},

\begin{assumption}
Observability matrix for $\mathbf{z}_s$, i.e., $\mathcal{O} \left( \mathbf{A}, \mathbf{H}_s \right) $ is rank deficient. 
\end{assumption}
\begin{subequations} \label{E:Ob1} \begin{align}
 \text{rank} \left[\begin{array}{c}
\mathbf{H}_s       \\
\mathbf{H}_s\mathbf{A}    \\
\mathbf{H}_s\mathbf{A} ^2     \\
\dots \\
\mathbf{H}_s\mathbf{A} ^{n-1}
\end{array}\right]  &= n_1 < n
\end{align} \end{subequations}

Using a suitable transformation matrix $\mathbf{\mathcal{T}}$, we can obtain a canonical form for such a system wherein the transformed state-space is partitioned into the observable and unobservable sub-spaces relative to the `suspicious measurements' \cite{kailath1980linear, Rugh}. The linearized dynamics of the noisy power systems in the transformed space can be expressed as 
\begin{align}
    \mathbf{\xbigbig}(k+1) &= \mathbf{\mathcal{A}}\,\mathbf{\xbigbig}(k)+ \mathbf{\mathcal{W}}(k)
\end{align}
with measurements 
\begin{align}
    \mathbf{z}(k) = \mathbf{\mathcal{H}}\,\mathbf{\xbigbig}(k) + \mathbf{\mathcal{V}}(k)\,,
\end{align} 
where 
\begin{align}
    \mathbf{\xbigbig}(k)^T=[(\mathbf{\xbigbig}_o^s)^T~ (\mathbf{\xbigbig}_u^s)^T]
\end{align} 
is the transformed state vector decomposed into the observable ($\mathbf{\xbigbig}_0^s$) and unobservable ($\mathbf{\xbigbig}_u^s$) parts with respect to $\mathbf{z}_s$, while the transformed system matrix 
\begin{align}
    \mathcal{\mathbf{A}}=\begin{bmatrix}\mathbf{\mathcal{A}}_{o} & \mathbf{0}\\\mathbf{\mathcal{A}}_{21}&\mathbf{\mathcal{A}}_{u}\end{bmatrix}
\end{align}
is also partitioned appropriately. We can use this transformed system description to distinguish between the two aforementioned scenarios.
Note that the following holds,
\begin{proposition} \label{T:111}
A change in $z_s$ will not affect ${\mathbf{\hat{\mathcal{X}}}}_{u}^s$. 
\end{proposition}
\noindent\begin{proof}
From the geometric interpretation of observablity \cite{kailath1980linear},
\begin{align}
    \mathbf{\mathcal{\hat{X}}}_{u}^s \in \mathcal{N}( \mathbf{\mathcal{O}} ( \mathbf{\mathcal{A}}, \mathcal{\mathbf{H}}_s))^T.
\end{align}
Since, the  states   ${\mathbf{{\mathcal{X}}}}_{u}^s$ lies in the null space of observablity matrix with respect to $\mathbf{z}_s$, a deviation in $\mathbf{z}_s$ will not affect their estimation  ${\mathbf{\hat{\mathcal{X}}}}_{u}^s$.
\end{proof}

Therefore, from Proposition\,\ref{T:111}, it follows that, in case of scenario-1, only the estimation of the observable  states $ \hat{\mathbf{\xbigbig}}_{o}^s$ are affected, manifested as a large deviation from the model based prediction $\mathbf{\mathcal{A}}_{o} \hat{\mathbf{\xbigbig}}_o^s(k+1|k)$, while the estimation of the unobservable states should remain unaffected. 

Let us define a metric $d(k+1)$.
\begin{definition}
The metric $d(k+1)$ is the $\mathcal{L}_{2}$-norm of the error between the a-posteriori and a-priori estimation of the unobservable states $\mathbf{\mathcal{\hat{X}}}_{u}^s$ at time instant $(k+1)$. Therefore, from \eqref{E:trans}, it follows:
\begin{align}
    &d(k+1)=\notag \\&\| \hat{\mathbf{\xbigbig}}_u^s(k+1|k+1) - \mathbf{\mathcal{A}}_{21} \hat{\mathbf{\xbigbig}}_o^s(k+1|k)- \mathbf{\mathcal{A}}_{u} \hat{\mathbf{\xbigbig}}_u^s(k+1|k)\|_2.
\end{align}
\end{definition}

\begin{figure}[thpb] 
  \centering
  \subfigure[Scenario-1: malicious measurements ($\mathbf{z}_s$)]{
  \includegraphics[scale=0.25]{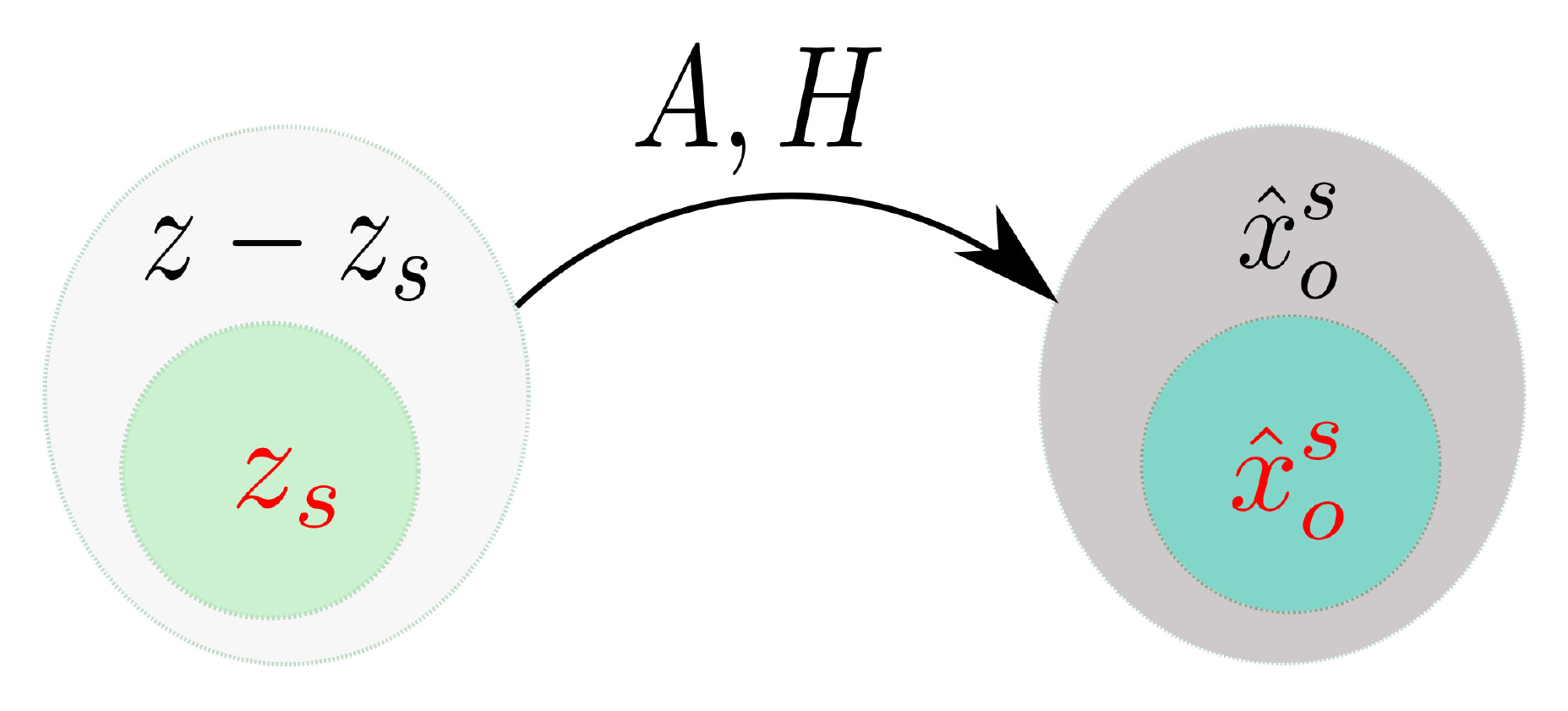}}
  \subfigure[Scenario-2: wrong model ($\mathbf{H},\,\mathbf{A}$).]{
  \includegraphics[scale=0.25]{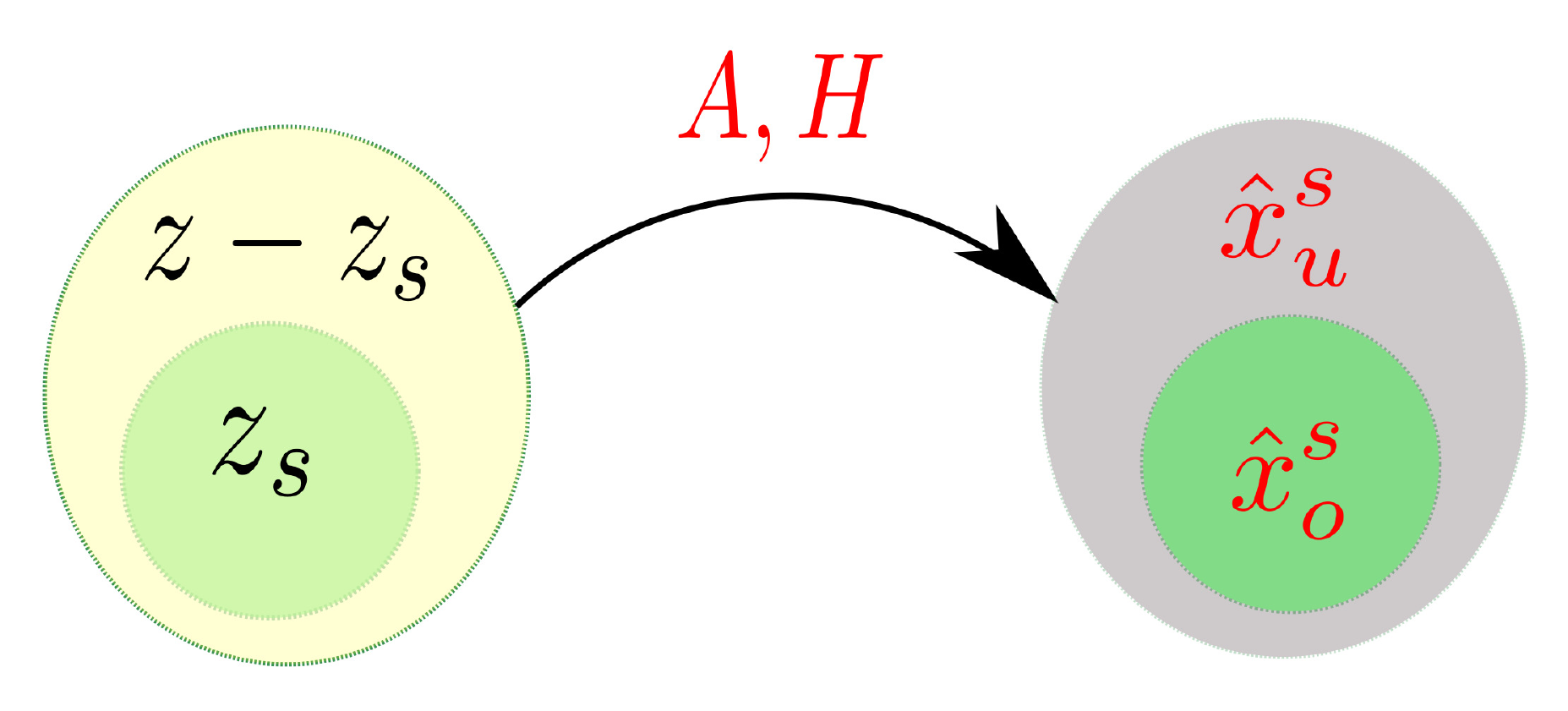}}
  \caption{Different scenarios of anomalies in estimation of system states. }
  \label{F:only}
\end{figure}

Therefore,  for scenario-1, $d(k+1)$ will be less.
Similarly, in case of scenario-2, both the observable and unobserbable states will be {affected} resulting into a larger $d(k+1)$.
Fig. \ref{F:only} shows, how the state estimations are affected for different scenarios. The affected variables are denoted in red.


 The following transform matrix $\mathbf{\mathcal{T}}$ will transform our state space model into observable and unobservable subsystems.
\begin{subequations}\begin{align}
\mathcal{\xbigbig}(k) & =\mathbf{\mathcal{T}} \xbigbig(k)\\
\left[\begin{array}{cccc}
{\mathcal{\xbigbig}}_{o}^s(k +1)      \\
{\mathcal{\xbigbig}}_{{u}}^s(k +1)    
\end{array}\right] &= \mathbf{\mathcal{A}}  \left[\begin{array}{cccc}
{\mathcal{\xbigbig}}_{o}^s (k )       \\
{\mathcal{\xbigbig}}_{{u}}^s (k )  
\end{array}\right]  + \mathbf{\mathcal{T}} \mathbf{w} (k ). \label{E:trans}
\end{align} \end{subequations}
Similarly, the measurement equation is:
\begin{subequations}\begin{align}
{\mathbf{z}}_s(k ) & = \mathbf{\mathcal{H}}_s \mathcal{\xbigbig} (k )+ \mathbf{\mathcal{T}} \mathbf{v}_s(k )\,,\\
\mathbf{\mathcal{H}}_s&=  [\mathbf{\mathcal{H}}_{s, o} \quad 0 ]
\end{align}\end{subequations}
One of the ways $\mathbf{\mathcal{T}}$ could be formed is \cite{kailath1980linear, Rugh}:
\begin{subequations}\begin{align}
\mathbf{\mathcal{T}} &= \left[\begin{array}{c} \mathbf{t}_1 \\ \mathbf{t}_2 \\ \vdots \\ \mathbf{t}_{n_1} \\ \mathbf{r}_1 \\ \mathbf{r}_2 \\ \vdots \\ \mathbf{r}_{n-n_1} \end{array}\right]  
\end{align} \end{subequations}
where, first $n_1$ rows are formed using $n_1$ independent columns from the observability matrix given by \eqref{E:Ob1}. Other $n-n_1$ rows are chosen arbitrarily so that $\mathbf{\mathcal{T}}$ is non-singular. 
The parameters at the transformed state space can be found to be,
\begin{subequations}\begin{align}
\mathbf{\mathcal{A}}&= \mathbf{\mathcal{T}}\mathbf{A}\mathbf{\mathcal{T}}^{-1}, \\  \mathbf{\mathcal{H}}_s &= \mathbf{H}_s \mathbf{\mathcal{T}}^{-1}, \\
\mathbf{\mathcal{Q}}&=\mathbf{\mathcal{T}}\mathbf{Q}\mathbf{\mathcal{T}}^{-1}.
\end{align} \end{subequations}
Now the system parameters corresponding to the observable subsystem:
\begin{subequations}\begin{align}
\mathbf{\mathcal{A}}_{o} & =\mathbf{\mathcal{A}} (1:n_1;1:n_1) \\
\mathbf{\mathcal{H}}_{s, o} & = \mathbf{\mathcal{H}}_s (1:m_s;1:n_1)
\end{align} \end{subequations}

Next, we transform the estimated state vector:
\begin{align}
\hat{\mathcal{\xbigbig}}=\mathbf{\mathcal{T}}\hat{\x}= 
\left[\begin{array}{c} \hat{\mathcal{\xbigbig}}_{o}^s \\ \dots \dots \\ \hat{\mathcal{\xbigbig}}_{u}^s    \end{array} \right]  
\end{align}
Now the transformed states unobservable to $\mathbf{z}_\mathbf{S}(k+1)$ are ${\hat{\mathcal{\xbigbig}}}_{{u}}^s(k+1|k+1)$, these estimated states are used to detect modeling error scenario using the following algorithm.
{\section{Algorithm}}
 For each time step $\left(k+1\right)$ when $c'(k+1) >TH_{\chi}$:
\begin{enumerate}
\item \textbf{Identify:} suspicious measurements $\mathbf{z}_s $, by $ \forall i, |\frac{z_i(k+1)-\mathbf{H}_i \hat{\mathbf{x}}(k+1|k)}{S_{i,i}(k+1)}| > TH_r $.
\item \textbf{Construct} $\O( \mathbf{A},\mathbf{z}_s)$
\item \textbf{If:} $rank ( \O(\mathbf{A},\mathbf{z}_s) ) = n_1< n$
\begin{enumerate}
\item \textbf{Construct:} $\mathbf{\mathcal{T}} \in \mathcal{R}^{n \times n}$ by choosing $n_1$ rows from $\O(A,\mathbf{z}_s)$ and arbitrarily choosing $n-n_1$ columns such that $\mathbf{\mathcal{T}}$ is full rank. 
\item \textbf{Compute:} $\mathbf{\hat{\mathcal{X}}}(k+1|k),\mathbf{\hat{\mathcal{X}}}(k+1|k+1),\mathbf{\mathcal{A}}$.
\item \textbf{Partition:} the transformed estimated states and system parameter and obtain $\mathbf{\hat{\mathcal{X}}_{u}}(k+1|k+1),\mathbf{\hat{\mathcal{X}}}_{u}(k+1|k), {\mathbf{\mathcal{A}}}_{21},{\mathbf{\mathcal{A}}}_{u}, \mathbf{\hat{\mathcal{X}}}_{{o}}(k+1|k)$.
\item \textbf{Compute:} $d(k+1)=||\mathbf{\hat{\mathcal{X}}_{u}}(k+1|k+1)- {\mathbf{\mathcal{A}}}_{21} \mathbf{\hat{\mathcal{X}}}_{{o}}(k+1|k) -  {\mathbf{\mathcal{A}}}_{u} \mathbf{\hat{\mathcal{X}}}_{u}(k+1|k)||_2$
\item \textbf{Determine:} the the threshold $TH_d$ for $d(k+1)$ is determined.
\item \textbf{If:} $d(k+1) > TH_d$, detect modeling error as the cause of anomaly (scenario-1).
\item \textbf{else:} detect malicious data (scenario-2).
\end{enumerate}
\item \textbf{else:} 
\begin{enumerate} 
\item \textbf{If:} $m_s > N_{critical}$, modeling error detected (scenario-1).
\item       \textbf{else:} Undecided.
\end{enumerate}
\item \textbf{end}
\end{enumerate}

In the above algorithm, there are four thresholds that has been used. Here, the description and the selection of each of these thresholds are presented.

\begin{enumerate}

    \item $\mathbf{{TH_{\chi}}:}$ The threshold for initial bad data detection's $\chi^2$-test. This threshold is detected using the inverse cumulative $\chi^2$ distribution corresponding to chosen confidence value $p$. 
    
    \item $ \mathbf{{TH_{r}}:}$ Threshold for individual normalized measurement residuals $r_i \quad \forall i \in \lbrace 1, 2, \dots, m \rbrace $. 
    
    \item $ \mathbf{{TH_{d}}:}$ This is the threshold value of the state squared sum of error between filtered and predicted estimation of unobservable states. A threshold can be chosen choosing proper confidence value corresponding to an inverse cumulative probability distribution.
    
    \item $ \mathbf{N_{critical}:}$ If most of the the measurements show large residuals, it could be more likely due to large modeling error, as it is unlikely that most of the meters will be attacked at the same time. This critical number of measurements $N_{critical}$ depends on the particular power system network and should be designed by the systems experts of the particular network. For our simulation, we assumed $N_{critical}=\left \lceil{\frac{m}{2}}\right \rceil$.
    
\end{enumerate}

\section{Results}
\begin{figure}[H]
 \centering
	\includegraphics[scale=0.75]{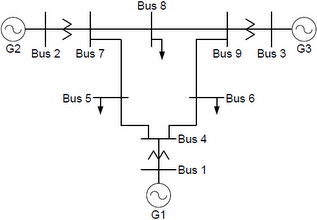}
 \caption{IEEE-9 bus test system. }
 \label{F:39}
 \end{figure}
 
\begin{figure}[thpb]
 \centering
	\includegraphics[scale=0.17]{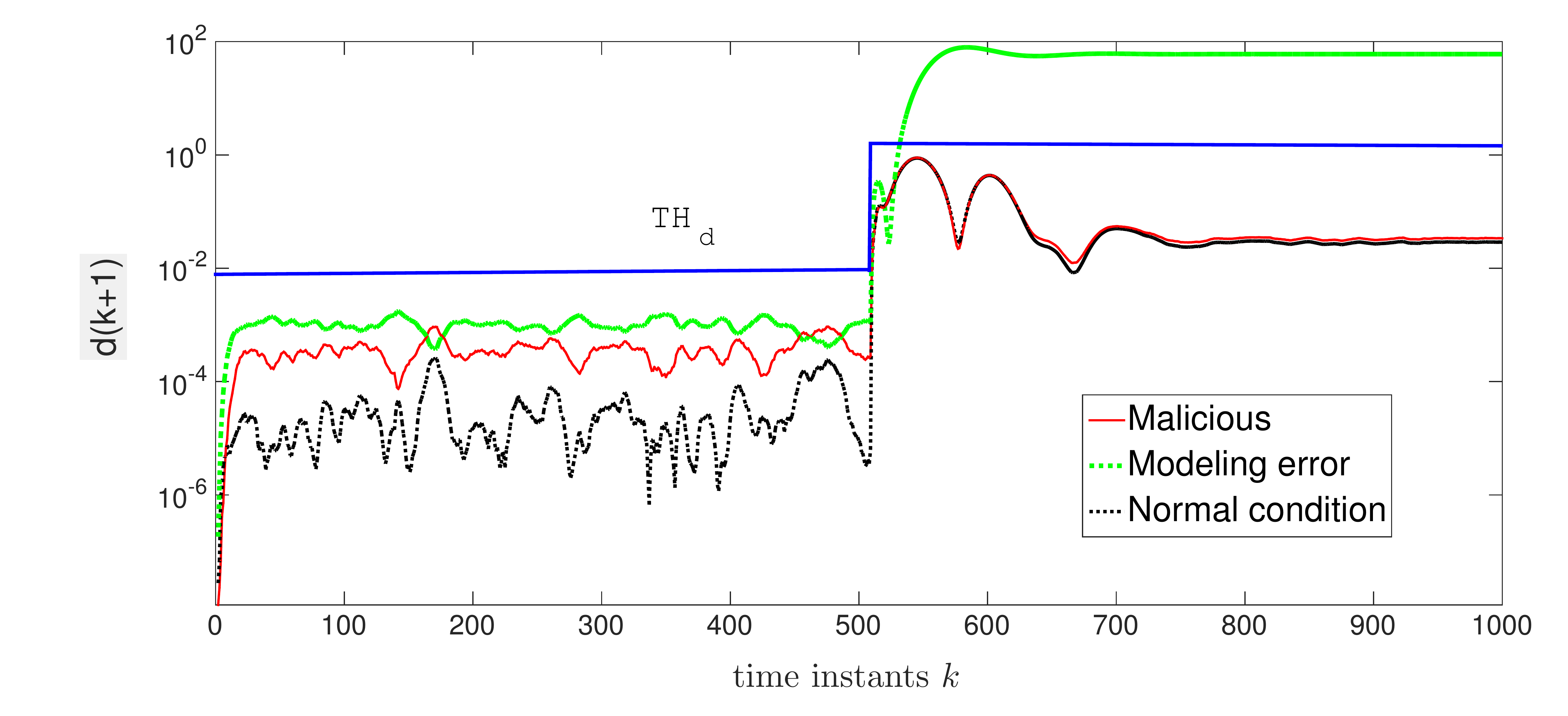}
 \caption{Comparison of $d(k+1)$ for malicious data, modeling error and normal condition. }
 \label{F:compare_anomaly}
 \end{figure}

{In this section IEEE-9 bus system is simulated. This is fourth order  system, having three generators. An L-L-L-G fault is applied at 5.1s-5.13s between the lines 5 and 7. PMU measurements are placed at the buses 1-8 getting voltage magnitude and phase angles. Then the simulation is repeated in three different scenarios:
\begin{itemize}
\item \textsc{ Modeling error (Scenario 2)}: Fault is applied but the model is not updated based on the new topology.
For this case, the normalized measurement residual and a threshold is plotted in Fig.\, \ref{F:mrr}. It can be observed that, after the fault is applied, the measurement residual, corresponding to the measurement-10 ($m_{10}$) exceeds the chosen threshold. 

\begin{figure}[tbh]
 \centering
	\includegraphics[scale=0.16]{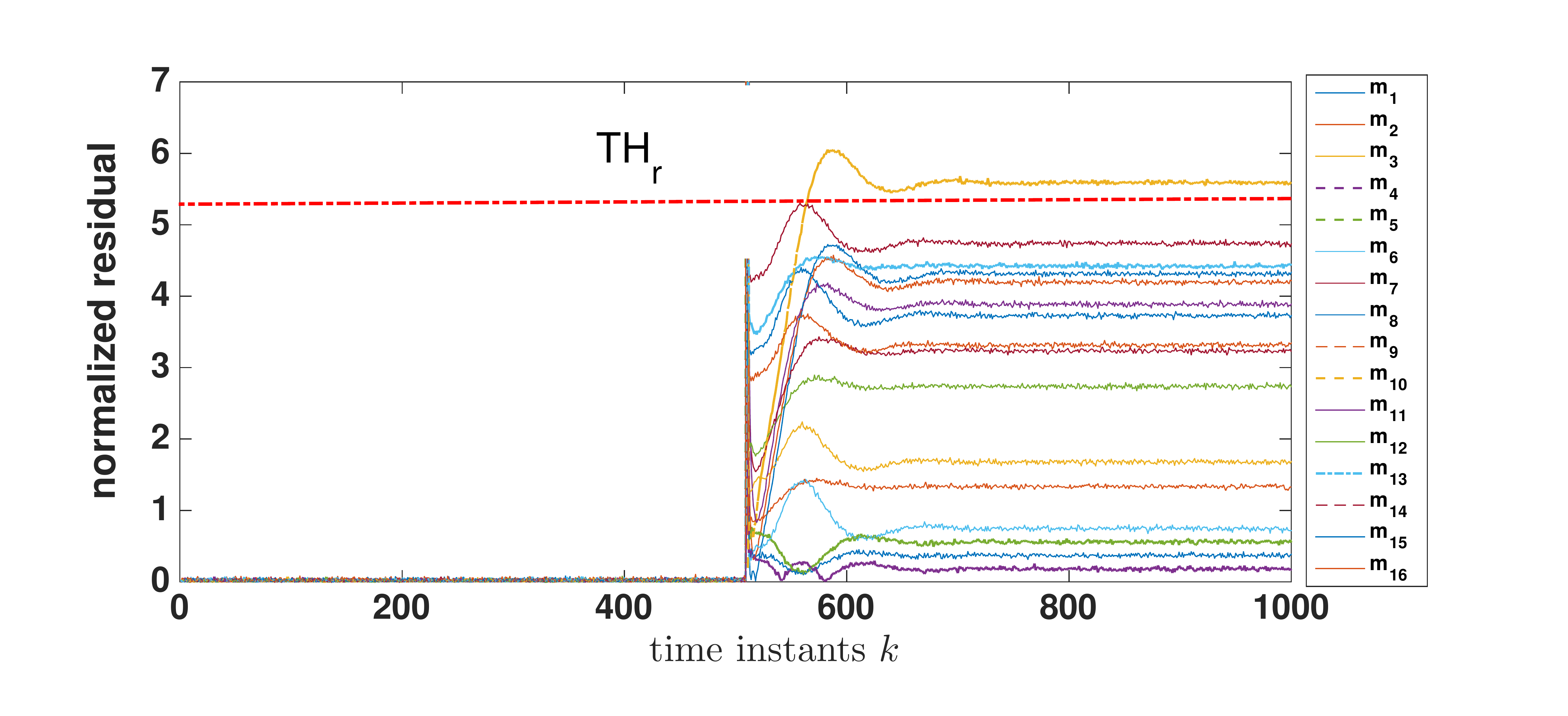}
 \caption{Normalized residuals corresponding to each measurements. }
 \label{F:mrr}
 \end{figure}

\item \textsc{Malicious data (Scenario 1) }: Fault is applied, model is updated based on the new topology and malicious data is injected into measurement-10.

\item \textsc{Normal condition}: Fault is applied and model is updated based on the new topology and no attack is injected.
\end{itemize}

A transformation matrix is identified to decompose the systems states into observable and unobservable parts, with respect to the measurement-10. Then, $d(k+1)$ is computed for the above three scenarios and plotted in Fig.\,\ref{F:compare_anomaly} (for all the scenarios). Clearly, the modeling error scenario is identified as $d(k+1)$ exceeds the threshold $TH_d$ during and after faults is applied.}

Please note that while the threshold $TH_r$ in Fig. \ref{F:mrr} is for a normalized quantity, in Fig. \ref{F:compare_anomaly} , the threshold $TH_d$ is not normalized, which changes with the change in operating condition.


 \section{Conclusion}
    
   The aim of the algorithm presented in this paper is to distinguish between the {cause} of the triggering of the bad-data detector:  1) Scenario-1: presence of malicious data and 2) Scenario-2: significant modeling error. This is important to diagnose the cause of the anomaly to  further reduce the false-alarms and also to improve the bad data detector performance. {This algorithm is applicable to any generic linear dynamical system.}

\IEEEtriggeratref{8}
\bibliographystyle{IEEEtran}
\bibliography{refelist}
\end{document}